\documentclass[a4paper,reqno]{amsart}

\usepackage[ngerman,english]{babel}
\usepackage{amsthm}
\usepackage{amsmath}
\usepackage{mathptmx}
\usepackage{amssymb}
\usepackage{txfonts}
\usepackage[usenames,dvipsnames]{xcolor}
\usepackage[colorlinks]{hyperref}
\usepackage{microtype}
\usepackage{enumitem}
\usepackage{pstool}
\usepackage{subfig}
\usepackage{cleveref}
\usepackage[utf8]{inputenc}
\usepackage{fullpage}
\usepackage[foot]{amsaddr}
\usepackage{mathrsfs}  

\theoremstyle{plain}
\newtheorem{theorem}{Theorem}[section]
\newtheorem{lemma}[theorem]{Lemma}
\newtheorem{corollary}[theorem]{Corollary}

\theoremstyle{definition}
\newtheorem{definition}[theorem]{Definition}
\newtheorem{remark}[theorem]{Remark}

\hypersetup{
	allcolors = Red,
	allbordercolors = Red,
	colorlinks = true}

\numberwithin{equation}{section}

\begin{document}

\title{On the density of intermediate $\beta$-shifts of finite type}

\author[B.\ Li]{Bing Li}
\address[Bing Li]{Department of Mathematics, South China University of Technology, Guangzhou, China}
\author[T.\ Sahlsten]{Tuomas Sahlsten}
\address[Tuomas Sahlsten]{School of Mathematics, The University of Manchester, Manchester, UK}
\author[T.\ Samuel]{Tony Samuel}
\address[Tony Samuel]{Mathematics Department, California Polytechnic State University, San Luis Obispo, CA, USA}
\author[W.\ Steiner]{Wolfgang Steiner}
\address[Wolfgang Steiner]{IRIF, CNRS, Universit\'e Paris Diderot - Paris 7, Paris, France}

 
\begin{abstract}
We determine the structure of the set of intermediate $\beta$-shifts of finite type.  Specifically, we show that this set is dense in the parameter space $\Delta = \{ (\beta, \alpha) \in \mathbb{R}^{2} \colon \beta \in (1, 2) \; \text{and} \; 0 \leq \alpha \leq 2 - \beta\}$.  This generalises the classical result of Parry from 1960 for greedy and (normalised) lazy  $\beta$-shifts.
\end{abstract}

\maketitle

\section{Introduction and statement of main results}

Since the pioneering work of R\'enyi \cite{R:1957} and Parry \cite{P:1960,P1964,P:1979}, intermediate $\beta$-shifts have been well-studied by many authors, and have been shown to have important connections with ergodic theory, fractal geometry and number theory.  We refer the reader to \cite{Komornik:2011,S:2003a} for survey articles concerning this topic.  Except for the results of \cite{LSS:2016,P:1960} and a few remarks in \cite{Par:1979,P:1979}, to the best of our knowledge, little is known concerning intermediate $\beta$-shifts of finite type, which is the focus of this article.

For $\beta > 1$ and $x \in [0,1/(\beta-1)]$, a word $(\omega_{n})_{n \in \mathbb{N}}$ with letters in the alphabet $\{ 0, 1\}$ is called a $\beta$-\textsl{expansion} of $x$ if
\begin{align*}
x = \sum_{k= 1}^{\infty} \omega_{k} \ \beta^{-k}.
\end{align*}
When $\beta$ is a natural number, Lebesgue almost all numbers $x$ have a unique $\beta$-expansion.  On the other hand, in \cite{S:2003b} it was shown that if $\beta$ is not a natural number, then, for almost all $x$, the cardinality of the set of $\beta$-expansions of $x$ is equal to the cardinality of the continuum.  In \cite{Baker:2015,BZWL:2016,SV:1998} the set of points which have countable number of $\beta$-expansions has also been studied.

Through iterating the maps $G_{\beta} \colon x \mapsto \beta x \bmod 1$ on $[0,1)$ and $L_{\beta} \colon x \mapsto \beta (x - 1) \bmod 1$ on $(0,1]$, for $\beta \in (1,2]$, one obtains subsets of $\{0, 1\}^{\mathbb{N}}$ known as the greedy and (normalised) lazy $\beta$-shifts, respectively, where each point $\omega^{+}$ of the greedy $\beta$-shift is a $\beta$-expansion, and corresponds to a unique point in $[0, 1]$, and each point $\omega^{-}$ of the lazy $\beta$-shift is a $\beta$-expansion, and corresponds to a unique point in $[(2-\beta)/(\beta - 1), 1/(\beta - 1)]$.  Note that, if $\omega^{+}$ and $\omega^{-}$ are $\beta$-expansions of the same point, then $\omega^{+}$ and $\omega^{-}$ do not necessarily have to be equal, see \cite{KL:1998}.

There are many ways, other than using the greedy and lazy $\beta$-shift, to generate a \mbox{$\beta$-expansion} of a positive real number. For instance the intermediate $\beta$-shifts $\Omega_{\beta, \alpha}$ which arise from the intermediate $\beta$-transformations $T^{\pm}_{\beta,\alpha} \colon [0,1] \circlearrowleft$, where $(\beta,\alpha) \in \Delta \coloneqq \{ (\beta, \alpha) \in \mathbb{R}^{2} \colon \beta \in (1, 2) \; \text{and} \; 0 \leq \alpha \leq 2 - \beta\}$, and where the maps $T_{\beta, \alpha}^{\pm}$ are defined as follows.  Letting $p = p_{\beta, \alpha} \coloneqq (1-\alpha)/\beta$ we set
\begin{align*}
T^{+}_{\beta, \alpha}(x) \coloneqq 
\begin{cases}
\beta x + \alpha \bmod 1 & \text{if} \; x \neq p,\\
0 & \text{if} \; x = p,
\end{cases}
\quad \text{and} \quad
T^{-}_{\beta, \alpha}(x) \coloneqq 
\begin{cases}
\beta x + \alpha \bmod 1 & \text{if} \; x \neq p,\\
1 & \text{if} \; x = p.
\end{cases}
\end{align*}
The maps $T_{\beta, \alpha}^{\pm}$ are equal everywhere except at the point $p$ and $T^{-}_{\beta, \alpha} (x) = 1 - T_{\beta, 2 - \beta - \alpha}^+(1 - x)$, for all $x \in [0,1]$. Notice, when $\alpha = 0$, the maps $G_{\beta}$ and $T^{+}_{\beta, \alpha}$ coincide on $[0,1)$, and when $\alpha = 2-\beta$, the maps $L_{\beta}$ and $T^{-}_{\beta, \alpha}$ coincide on $(0,1]$.  Intermediate $\beta$-transformations are sometimes called linear Lorenz maps and arise naturally from Poincar\'e maps of the geometric model of Lorenz differential equations, see for instance \cite{EO:1994,L:1963,V:2003,W:1980}.  Here, let us also make the observation that, for all $(\beta, \alpha) \in \Delta$, the symbolic space $\Omega_{\beta, \alpha}$ is always a subshift, meaning that it is invariant under the left shift map.  As an aside, we remark that another way to generate a $\beta$-expansion of a point is to use random $\beta$-transformations, see for instance \cite{DajKal:2007,DajKra:2003,DajVries:2005}.

Subshifts are to dynamical systems what shapes like polygons and curves are to geometry.  Subshifts which can be described by a finite set of forbidden words are called subshifts of finite type and play an essential role in the study of dynamical systems --- indeed, this class of subshifts play a similar role in dynamics as triangles play in geometry.  Their study has provided solutions to important practical problems, such as finding efficient coding schemes to store data \cite{LM:1995} and recently in the analysis of electroencephalography (EEG) data \cite{Stolz:2017}.

One reason why subshifts of finite type are so useful is that they have a simple representation using a finite directed graph.  Questions about the subshift can then often be phrased as questions about the graph's adjacency matrix.  Results from linear algebra help us to take this matrix apart and find solutions.  Moreover, in the case of greedy $\beta$-shifts (that is when $\alpha = 0$), to compute the multifractal analysis for Birkhoff averages, one first computes the result for greedy $\beta$-shifts of finite type, and then one uses an approximation argument to determine the result for a general greedy $\beta$-shift, see \cite{LL:2016} and references therein.

Given $(\beta, \alpha) \in \Delta$, the $\beta$-expansions of the point $p$ given by $T_{\beta, \alpha}^{\pm}$ are called the kneading invariants of $\Omega_{\beta, \alpha}$.  In \cite{BSV:2014,H:1979,HS:1990} it is shown that the kneading invariants completely determine $\Omega_{\beta, \alpha}$ and in \cite{LSS:2016} it is shown that, provided $\alpha \not\in \{ 0, 2-\beta \}$, then an intermediate $\beta$-shift is a subshift of finite type if and only if the kneading invariants are periodic; in contrast, the greedy and lazy $\beta$-shifts (that is when $\alpha = 0$ and $\alpha = 2 -\beta$, respectively) are subshifts of finite type if and only if one of the kneading invariants is periodic, see \cite{P:1960}.  Our contribution in this article to this story is to show the following result, which generalises the classical density result of Parry \cite{P:1960} for greedy and (normalised) lazy $\beta$-shifts.

\begin{theorem}\label{thm:density}
The set of $(\beta, \alpha)$ belonging to $\Delta$ for which $\Omega_{\beta, \alpha}$ is a subshift of finite type is dense in $\Delta$.
\end{theorem}

A subshift is \textsl{sofic} if it is a factor of a subshift of finite type.  \Cref{thm:density} implies that the set of $(\beta, \alpha) \in \Delta$ for which the intermediate $\beta$-shift is sofic is dense in $\Delta$, since any subshift of finite type is a sofic shift.  If one considers the dynamical property of topologically transitivity, then the structure of the set of $(\beta, \alpha)$ in $\Delta$ such that $\Omega_{\beta, \alpha}$ is topologically transitive, with respect to the left shift map, is very different to the set of $(\beta, \alpha)$ belonging to $\Delta$ for which $\Omega_{\beta, \alpha}$ is a subshift of finite type.  In fact the former of these two sets is far from being dense in $\Delta$, see \cite{G:1990,Par:1979}.

Combining the results of \cite{BHV:2011} and \cite{LSS:2016} it can be shown that if $\beta$ is a transcendental number, then the set of $\alpha$ for which the intermediate $\beta$-shift $\Omega_{\beta, \alpha}$ is a subshift of finite type is empty.  Indeed, for $\Omega_{\beta, \alpha}$ to be a subshift of finite type, we require $\beta \in (1, 2)$ to be a maximal root of a polynomial with coefficients in $\{ -1, 0, 1 \}$.  Moreover, as the entropy of a subshift of finite type is the logarithm of the largest eigenvalue of its adjacency matrix, which is a nonnegative integral matrix, we have the following.  If an intermediate $\beta$-shift is of finite type, then there exists an $n \in \mathbb{N}$ such that $\beta$ is the positive $n^{\textup{th}}$-root of a Perron number \cite{Lind:84}. This leads to the following very natural question.  If $\beta \in (1, 2)$ is a positive $n^{\textup{th}}$-root of a Perron number, for some $n \in \mathbb{N}$, then is the set of $\alpha$ for which $\Omega_{\beta, \alpha}$ is a subshift of finite type dense in  $[0, 2 - \beta]$?

\section{Preliminaries}

\subsection{Subshifts}

We equip the space $\{0,1\}^\mathbb{N}$ of infinite sequences with the topology induced by the word metric $\mathscr{D} \colon \{0,1\}^\mathbb{N} \times \{0,1\}^\mathbb{N} \to \mathbb{R}$ which is given by
\begin{align*}
\mathscr{D}(\omega, \nu) \coloneqq
\begin{cases}
0 & \text{if} \; \omega = \nu,\\
2^{- \lvert\omega \wedge \nu\rvert + 1} & \text{otherwise}.
\end{cases}
\end{align*}
Here $\rvert \omega \wedge \nu \lvert \coloneqq \min \, \{ \, n \in \mathbb{N} \colon \omega_{n} \neq \nu_n \}$, for all $\omega = (\omega_{1}, \omega_{2}, \dots) , \nu = ( \nu_{1} \, \nu_{2}, \dots) \in \{0, 1\}^{\mathbb{N}}$ with $\omega \neq \nu$. Note that the topology induced by $\mathscr{D}$ on $\{ 0, 1\}^{\mathbb{N}}$ coincides with the product topology on $\{ 0, 1\}^{\mathbb{N}}$.  We let $\sigma \colon \{ 0, 1 \}^{\mathbb{N}} \circlearrowleft$ denote the \textsl{left-shift map} which is defined by $\sigma(\omega_{1}, \omega_{2}, \dots) \coloneqq (\omega_{2}, \omega_{3}, \dots)$.  A \textsl{subshift} is any closed set $\Omega \subseteq \{0,1\}^\mathbb{N}$ such that $\sigma(\Omega) \subseteq \Omega$.  Given a subshift $\Omega$ and $n \in \mathbb{N}$ we set
\begin{align*}
\Omega\lvert_{n} \coloneqq \left\{ (\omega_{1}, \dots, \omega_{n}) \in \{ 0, 1\}^{n} \colon \,\text{there exists} \; (\xi_{1}, \xi_{2}, \dots ) \in \Omega \; \text{with} \; (\xi_{1}, \dots, \xi_{n}) = (\omega_{1}, \dots, \omega_{n}) \right\}
\end{align*}
and write $\Omega^{*} \coloneqq \bigcup_{n = 1}^{\infty} \Omega\lvert_{n}$ for the collection of all finite words.  We denote by $\lvert \Omega\vert_{n} \rvert$ the cardinality of $\Omega\vert_{n}$.  A subshift $\Omega$ is a called a \textsl{subshift of finite type} if there exists $M \in \mathbb{N}$ such that, $(\omega_{n - M + 1}, \dots, \omega_{n}, \xi_{1}, \dots, \xi_{m}) \in \Omega^{*}$, for all $(\omega_{1}, \dots, \omega_{n})$, $(\xi_{1}, \dots, \xi_{m}) \in \Omega^{*}$ with $n, m \in \mathbb{N}$ and $n \geq M$, if and only if $(\omega_{1}, \dots, \omega_{n}, \xi_{1}, \dots, \xi_{m}) \in \Omega^{*}$.

The following result gives an equivalent condition for when a subshift is of finite type.  For this we require the following notation.   For $n \in \mathbb{N}$ and $\omega = (\omega_{1}, \omega_{2}, \dots) \in \{ 0, 1\}^{\mathbb{N}}$, we set $\omega\lvert_{n} \coloneqq (\omega_{1}, \dots, \omega_{n})$ and, for $\xi \in \{ 0, 1 \}^{*}$, we let $\lvert \xi \rvert$ denote the length of $\xi$.  

\begin{theorem}[{\cite{LM:1995}}]\label{thm:equivalent-def-SFT}
A shift space $\Omega$ is a subshift of finite type if and only if there exists a finite set $F \subset \Omega^{*}$ such that, if $\xi \in F$, then $\omega\lvert_{\lvert \xi \rvert} \neq \xi$, for all $\omega \in \Omega$.  The set $F$ is often called the set of forbidden words.
\end{theorem}

For $n, m \in \mathbb{N}$ and $\nu = (\nu_{1}, \dots, \nu_{n}),\, \xi = (\xi_{1}, \dots, \xi_{m}) \in \{ 0, 1\}^{*}$, set $(\nu, \xi) \coloneqq (\nu_{1}, \dots, \nu_{n}, \xi_{1}, \dots, \xi_{m})$; we use the same notation when $\xi \in \{ 0, 1\}^{\mathbb{N}}$.

An infinite word $\omega = (\omega_{1}, \omega_{2}, \dots) \in \{0, 1\}^{\mathbb{N}}$ is called \textsl{periodic} with \textsl{period} $n \in \mathbb{N}$ if and only if, for all $m \in \mathbb{N}$, we have $(\omega_{1}, \dots, \omega_{n}) = (\omega_{(m - 1)n + 1}, \dots, \omega_{m n})$; in which case we write $\omega = (\overline{\omega_{1}, \dots, \omega_{n}})$.  Similarly, $\omega = (\omega_{1}, \omega_{2}, \dots) \in \{0, 1\}^{\mathbb{N}}$ is called \textsl{eventually periodic} with \textsl{period} $n \in \mathbb{N}$ if and only if there exists $k \in \mathbb{N}$ such that, for all $m \in \mathbb{N}$, we have $(\omega_{k+1}, \dots, \omega_{k+n}) = (\omega_{k+(m - 1)n + 1}, \dots, \omega_{k+ m n})$; in which case we write $\omega = (\omega_{1}, \dots, \omega_{k}, \overline{\omega_{k+1}, \dots, \omega_{k+n}})$.

\subsection{Intermediate $\beta$-shifts and expansions}

We now give the formal definition of an intermediate $\beta$-shift. Throughout this section we let $(\beta, \alpha) \in \Delta$ be fixed and let $p = p_{\beta,\alpha} \coloneqq (1 - \alpha)/\beta$.

The \textsl{$T^{\pm}_{\beta, \alpha}$-expansion} $\tau_{\beta, \alpha}^{\pm}(x)$ of $x \in [0, 1]$ is defined to be the infinite word $(\omega^{\pm}_{1}, \omega^{\pm}_{2}, \dots, ) \in \{ 0, 1\}^{\mathbb{N}}$, where, for $n \in \mathbb{N}$,
\begin{align*}
\omega^{+}_{n} \coloneqq \begin{cases}
0 & \quad \text{if } (T^{+}_{\beta,\alpha})^{n-1}(x) < p,\\
1 & \quad \text{otherwise,}
\end{cases}
\quad \text{and} \quad
\omega^{-}_{n} \coloneqq \begin{cases}
0 & \quad \text{if } \, (T^{-}_{\beta,\alpha})^{n-1}(x) \leq p,\\
1 & \quad \text{otherwise.}
\end{cases}
\end{align*}
We will denote the images of the unit interval under $\tau_{\beta, \alpha}^{\pm}$ by $\Omega^{\pm}_{\beta, \alpha}$, respectively, and set $\Omega_{\beta, \alpha} \coloneqq \Omega_{\beta, \alpha}^{+} \cup \Omega_{\beta, \alpha}^{-}$.  The \textsl{upper} and \textsl{lower kneading invariants} of $\Omega_{\beta,\alpha}$ are defined to be the infinite words $\tau^{\pm}_{\beta, \alpha}(p)$, respectively.

\begin{remark}\label{rem:00..11..}
Let $\omega^{\pm} = (\omega_{1}^{\pm}, \omega_{2}^{\pm}, \dots, )$ denote the infinite words $\tau_{\beta, \alpha}^{\pm}(p)$, respectively.  By definition, $\omega^{-}_{1} = \omega^{+}_{2} = 0$ and $\omega^{+}_{1} = \omega^{-}_{2} = 1$.  One can also show that $(\omega_{k}^{+}, \omega_{k+1}^{+}, \dots, ) = (0, 0, 0, \dots )$ if and only if $\alpha = 0$ and $k \ge 2$; and $(\omega_{k}^{-}, \omega_{k+1}^{-}, \dots, ) = (1, 1, 1, \dots )$ if and only if $\alpha = 2-\beta$ and $k \ge 2$.  Indeed, $\tau_{\beta, \alpha}^{+}(0) = \sigma(\tau_{\beta, \alpha}^{+}(p))$ and $\tau_{\beta, \alpha}^{-}(1) = \sigma(\tau_{\beta, \alpha}^{-}(p))$.
\end{remark}

The connection between intermediate $\beta$-transformations and the $\beta$-expansions of real numbers is given via the $T_{\beta, \alpha}^{\pm}$-expansions of a point and the projection $\pi_{\beta, \alpha} \colon \{ 0, 1 \}^{\mathbb{N}} \to [0, 1]$ defined by
\begin{align*}
\pi_{\beta, \alpha}(\omega_{1}, \omega_{2}, \dots) \coloneqq \frac{\alpha}{1 - \beta} + \sum_{k = 1}^{\infty} \frac{\omega_{k}}{\beta^k}.
\end{align*}
We note that $\pi_{\beta, \alpha}$ is linked to the iterated function system $( [0, 1]; f_{0} \colon x \mapsto x/\beta, f_{1} \colon x \mapsto (x+1)/\beta)$ via the equality
\begin{align*}
\pi_{\beta, \alpha}(\omega_{1}, \omega_{2}, \dots) = \alpha / (1-\beta) + \lim_{n \to \infty} f_{\omega_{1}} \circ \dots \circ f_{\omega_{n}}([0, 1]).
\end{align*}
We refer the reader to \cite{F:1990,F:1997} for further details on iterated function systems.  An important property $\pi_{\beta, \alpha}$ is that the following diagram commutes.
\begin{align*}
\begin{array}
[c]{ccc}
\Omega_{\beta, \alpha}^{\pm} & \overset{\sigma}{\longrightarrow} & \Omega_{\beta, \alpha}^{\pm}\\
& & \\
\pi_{\beta, \alpha} \downarrow \uparrow \tau^{\pm}_{\beta, \alpha} &  & \pi_{\beta, \alpha} \downarrow \uparrow \tau^{\pm}_{\beta, \alpha} \\ & & \\
 \lbrack 0,1 \rbrack & \underset{T_{\beta, \alpha}^{\pm}} {\longrightarrow} & \lbrack 0,1\rbrack
\end{array}
\end{align*}
This result is readily verifiable from the definitions of the maps involved and a sketch of a proof can be found in \cite{BHV:2011}.  From this, one can deduce that, for $x \in [\alpha/(\beta - 1), 1 + \alpha/(\beta-1)]$, the words $\tau^{\pm}_{\beta, \alpha}(x - \alpha /(\beta - 1))$ are $\beta$-expansions of~$x$.

\subsection{Intermediate $\beta$-shifts of finite type}

An equivalent condition, in terms of the upper and lower kneading invariants, for the finite type property of greedy and the lazy $\beta$-shifts is given in \cite{P:1960} and reads as follows.

\begin{theorem}[{\cite{P:1960}}]\label{thm:Parry1960}
For $\beta \in (1, 2)$, we have that
\begin{enumerate}[leftmargin=2.5em]
\item\label{thmA:2} the greedy $\beta$-shift (that is when $\alpha = 0$) is a subshift of finite type if and only if $\tau^{-}_{\beta, 0}(p)$ is periodic, and
\item\label{thmA:3} the lazy $\beta$-shift (that is when $\alpha = 2 - \beta$) is a subshift of finite type if and only if $\tau^{+}_{\beta, 2-\beta}(p)$ is periodic.
\end{enumerate}
\end{theorem}

An analogous result holds for intermediate $\beta$-shifts.

\begin{theorem}[{\cite{LSS:2016}}]\label{thm:ESFTP}
Let $\beta \in (1, 2)$ and $\alpha \in (0, 2 - \beta)$.  The intermediate $\beta$-shift $\Omega_{\beta, \alpha}$ is a subshift of finite type if and only if both $\tau^{\pm}_{\beta, \alpha}(p)$ are periodic.
\end{theorem}

A necessary and sufficient condition, in terms of the upper and lower kneading invariants, for the property of an intermediate $\beta$-shift to be a sofic shift can be found in \cite{KS:2012} and reads as follows.  Let $(\beta, \alpha) \in \Delta$.  The intermediate $\beta$-shift $\Omega_{\beta, \alpha}$ is sofic if and only if both $\tau^{\pm}_{\beta, \alpha}(p)$ are eventually periodic.

\subsection{Structure of intermediate $\beta$-shifts}

The following results on the structure of $\Omega_{\beta, \alpha}^{\pm}$ will play a crucial role in our proof of \Cref{thm:density}.

\begin{theorem}[{\cite{AM:1996,BHV:2011,BM:2011,HS:1990,KS:2012}}]\label{thm:Structure}
For $(\beta, \alpha) \in \Delta$, the spaces $\Omega_{\beta, \alpha}^{\pm}$ are completely determined by upper and lower kneading invariants of $\Omega_{\beta, \alpha}$;  indeed, we have
\begin{align*}
\Omega_{\beta, \alpha}^{+} &= \left\{ \omega \in \{ 0, 1\}^{\mathbb{N}} \colon \tau_{\beta, \alpha}^{+}(0) \preceq \sigma^{n}(\omega) \prec \tau_{\beta, \alpha}^{-}(p) \; \textup{or} \; \tau_{\beta, \alpha}^{+}(p) \preceq \sigma^{n}(\omega) \preceq \tau_{\beta, \alpha}^{-}(1) \; \text{for all} \; n \in \mathbb{N}_{0} \right\},\\
\Omega_{\beta, \alpha}^{-} &= \left\{ \omega \in \{ 0, 1\}^{\mathbb{N}} \colon \tau_{\beta, \alpha}^{+}(0) \preceq  \sigma^{n}(\omega) \preceq \tau_{\beta, \alpha}^{-}(p) \; \textup{or} \; \tau_{\beta, \alpha}^{+}(p) \prec \sigma^{n}(\omega) \preceq \tau_{\beta, \alpha}^{-}(1) \; \text{for all} \; n \in \mathbb{N}_{0} \right\}.
\end{align*}
Here, $\prec$, $\preceq$, $\succ$ and $\succeq$ denote the lexicographic orderings on $\{ 0 ,1\}^{\mathbb{N}}$.  Moreover, $\Omega_{\beta, \alpha} = \Omega_{\beta, \alpha}^{+} \cup \Omega_{\beta, \alpha}^{-}$ is closed with respect to the metric $\mathscr{D}$ and hence is a subshift.
\end{theorem}

A natural question stemming from this result is the following.  When are two elements of $\{ 0, 1 \}^{\mathbb{N}}$ kneading invariants of an intermediate $\beta$-transformation?  This was answered in \cite{BSV:2014,HS:1990}.  To state the result we require the following.  Given $\omega, \nu \in \{ 0, 1\}^{\mathbb{N}}$ with $\omega \prec \nu$  we define the sets $\Omega^{\pm}(\omega, \nu)$ by
\begin{align*}
\Omega^{+}(\omega, \nu) &\coloneqq \{ \xi \in \{ 0, 1 \}^{\mathbb{N}} \colon \sigma(\nu) \preceq \sigma^{n}(\xi) \prec \omega \; \text{or} \; \nu \preceq \sigma^{n}(\xi) \preceq \sigma(\omega) \; \text{for all} \; n \in \mathbb{N}_{0}\},\\
\Omega^{-}(\omega, \nu) &\coloneqq \{ \xi \in \{ 0, 1 \}^{\mathbb{N}} \colon \sigma(\nu) \preceq \sigma^{n}(\xi) \preceq \omega \; \text{or} \; \nu \prec \sigma^{n}(\xi) \preceq \sigma(\omega) \; \text{for all} \; n \in \mathbb{N}_{0} \}.
\end{align*}

\begin{definition}\label{defn:admissible}
A pair $(\omega, \nu)$ of infinite words is said to be \textsl{admissible} if the following four conditions are satisfied.
\begin{enumerate}[leftmargin=2.5em]
\item\label{enumi1:defn_admissible} $\omega\lvert_{1} = 0$ and $\nu\lvert_{1} = 1$
\item\label{enumi2:defn_admissible} $\omega \in \Omega^{-}(\omega, \nu)$ and $\nu \in \Omega^{+}(\omega, \nu)$
\item\label{enumi3:defn_admissible} $\displaystyle{\lim_{n \to \infty} n^{-1} \ln \left( \lvert \Omega(\omega, \nu)\lvert_{n} \rvert \right) > 0}$
\item\label{enumi4:defn_admissible} If $\omega, \nu \in \{ \xi, \zeta \}^{\mathbb{N}}$ for two finite words $\xi, \zeta$ in the alphabet $\{ 0, 1\}$ with length greater than or equal to three, such that $\xi\lvert_{2} = (0,1)$, $\zeta\lvert_{2} = (1,0)$, $(\overline{\xi}) \in \Omega^{-}((\overline{\xi}), (\overline{\zeta}))$ and $(\overline{\zeta}) \in \Omega^{+}((\overline{\xi}), (\overline{\zeta}))$, then $\omega = (\overline{\xi})$ and $\nu = (\overline{\zeta})$.
\end{enumerate}
If in addition $\omega$ and $\nu$ are periodic words, then we call the pair $(\omega, \nu)$ \textsl{periodically admissible}.

\end{definition}

\begin{remark}\label{rmk:subaddiitve}
The limit in \Cref{defn:admissible} \eqref{enumi3:defn_admissible} exists for the following reason.  For a given pair $(\omega, \nu)$ of infinite words and $m, n \in \mathbb{N}$, we have that $\lvert \Omega(\omega, \nu)\vert_{m + n} \rvert \leq \lvert \Omega(\omega, \nu)\vert_{m} \rvert \ \lvert \Omega(\omega, \nu)\vert_{n} \rvert$. Therefore, the sequence $( \ln(\lvert \Omega(\omega, \nu)\vert_{n} \rvert)_{n \in \mathbb{N}}$ is sub-additive and hence $\displaystyle \lim_{n \to \infty} n^{-1} \ln \left( \lvert \Omega(\omega, \nu)\lvert_{n} \rvert \right) = \inf \{  n^{-1} \ln \left( \lvert \Omega(\omega, \nu)\lvert_{n} \rvert \right) \colon n \in \mathbb{N} \}$.
\end{remark}

With this at hand we may answer the above question.

\begin{theorem}[\cite{BSV:2014,HS:1990}]\label{thm:BSV2014}
Two infinite words $\omega, \nu \in \{ 0, 1\}^{\mathbb{N}}$ are kneading invariants for an intermediate $\beta$-shift if and only if $(\omega, \nu)$ is an admissible pair.
\end{theorem}

\subsection{Periodic kneading invariants}

\begin{corollary}\label{cor:Lss2016_BSV2014}
Two infinite words $\omega, \nu \in \{ 0, 1\}^{\mathbb{N}}$ are kneading invariants for an intermediate $\beta$-shift of finite type if and only if $(\omega, \nu)$ is a periodically admissible pair.
\end{corollary}

The following lemma, which is interesting in its own right and which we will require in our proof of \Cref{thm:density}, shows that \Cref{defn:admissible} \eqref{enumi4:defn_admissible} is violated when \Cref{defn:admissible} \eqref{enumi1:defn_admissible}, \eqref{enumi2:defn_admissible} and \eqref{enumi3:defn_admissible} hold but where $\omega$ and $\nu$ are not kneading invariants of an intermediate $\beta$-shift.

\begin{lemma}[{\cite{BSV:2014}}]\label{lem:BSV:2014}
Let $\omega = (\omega_{1}, \omega_{2}, \dots ), \nu = ( \nu_{1}, \nu_{2}, \dots ) \in \{ 0, 1\}^{\mathbb{N}}$ be such that \Cref{defn:admissible} \eqref{enumi1:defn_admissible}, \eqref{enumi2:defn_admissible} and \eqref{enumi3:defn_admissible} are satisfied.  Let $\beta$ denote the exponential of the value given in  \Cref{defn:admissible} \eqref{enumi3:defn_admissible} and let
\begin{align}\label{eq:BSV:2014_alpha}
\alpha = 1 - \beta + \beta (\beta - 1) \sum_{k = 1}^{\infty} \frac{\omega_{k}}{\beta^{k}},
\end{align}
namely the solution of $\pi_{\beta, \alpha}(\omega)  = p$, with $p = (1 - \alpha)/\beta$. 
If $\omega \neq \tau^{-}_{\beta, \alpha}(p)$ or $\nu \neq \tau^{+}_{\beta, \alpha}(p)$, then there exist $u, v \in \{ 0, 1\}^{*}$ such that $\omega, \nu \in \{ u, v \}^{\mathbb{N}}$ and $((\overline{u}), (\overline{v}))$ satisfies \Cref{defn:admissible} \eqref{enumi1:defn_admissible}, \eqref{enumi2:defn_admissible} and \eqref{enumi3:defn_admissible} with $\omega \neq (\overline{u})$ or $\nu \neq (\overline{v})$.  Moreover, $(\overline{u}) = \tau^{-}_{\beta, \alpha}(p)$ and $(\overline{v}) = \tau^{+}_{\beta, \alpha}(p)$.

\end{lemma}

\section{Proof of \Cref{thm:density}}

Fix $(\beta, \alpha) \in \Delta$ with $\alpha \not\in \{ 0, 2-\beta \}$.  For ease of notation let $\nu = (\nu_{1}, \nu_{2}, \dots ) = \tau_{\beta, \alpha}^{+}(p_{\beta, \alpha})$ and $\omega = (\omega_{1}, \omega_{2}, \dots ) =\tau_{\beta, \alpha}^{-}(p_{\beta, \alpha})$.  Assume that at least one of $\omega$ and $\nu$ are not periodic.  Our aim is to construct a pair of periodically admissible words $(\mu, \eta)$, arbitrarily close to the pair $(\omega, \nu)$. Then, by \Cref{cor:Lss2016_BSV2014} and \Cref{lem:BSV:2014}, there exists a unique point $(b, a) \in \Delta$ such that $\mu = \tau_{b, a}^{-}(p_{b, a})$, $\eta = \tau_{b, a}^{+}(p_{b, a})$, $\Omega_{b, a}$ is a subshift of finite type.  To complete the proof we show that the Euclidean distance between $(b, a)$ and $(\beta, \alpha)$ is arbitrarily small.

In what follows, we use that, by \Cref{thm:Structure}, for $m \in \mathbb{N}$, if $\omega_{m+1} = 0$ then $\sigma^m(\omega) \preceq \omega$, and if $\nu_{m+1} = 1$, then $\sigma^m(\nu) \succeq \nu$.
Suppose $\omega$ is not periodic, otherwise, set $\omega' = \omega$ and move to the construction of $\nu'$ in the following paragraph.  By \Cref{rem:00..11..}, we may choose an integer $n \geq 3$ with $\omega_{n} = 0$. 
Let $j \ge 1$ be minimal such that $(\omega_{j+1}, \dots, \omega_{n}) = (\omega_{1}, \dots, \omega_{n-j})$; note this equation holds for $j = n-1$. Since $\omega$ is not periodic, we have $\sigma^{j}(\omega) \prec \omega$. 
Therefore, there exists a minimal integer $k \geq n$ with $\omega_{k+1} = 0$ and $\omega_{k-j+1} = 1$, in which case,
\begin{align*}
(\omega_{j+1},\dots,\omega_{n-1}, \omega_n,\dots,\omega_{k}) = (\omega_1,\dots,\omega_{n-j-1},\omega_{n-j},\dots,\omega_{k-j}).
\end{align*}
Let $\omega' = (\overline{\omega_{1}, \dots, \omega_{k}})$.  
We claim that $(\omega', \nu)$ is a pair of infinite words satisfying \Cref{defn:admissible} \eqref{enumi1:defn_admissible} and \eqref{enumi2:defn_admissible}.  
Since $\sigma^{ik}(\omega)\vert_{k+1} \preceq \omega\vert_{k+1} = \omega'\vert_{k+1} = \sigma^{ik}(\omega')\vert_{k+1}$ for all $i \in \mathbb{N}$ with $\omega_{ik+1} = 0$, we have $\omega \preceq \omega'$.
As $\omega \neq \omega'$, it follows that $\omega \prec \omega'$. Combing this with the fact that $\omega\vert_{k+1} = \omega'\vert_{k+1}$, we conclude that $\sigma^m(\omega) \prec \sigma^m(\omega')$, for all $m \in \{0,\ldots,k+1\}$.
By the admissibility of $(\omega, \nu)$, we have $\nu \in \Omega^{+}(\omega, \nu)$ and, hence, $\nu \in \Omega^{+}(\omega', \nu)$.
For proving $\omega' \in \Omega^{-}(\omega', \nu)$, by the periodicity of~$\omega'$ and the definition of $\Omega^{-}(\omega', \nu)$, it is sufficient to check that $\sigma(\nu) \preceq \sigma^{m}(\omega') \preceq \omega'$ or  $\nu \prec \sigma^{m}(\omega') \preceq \sigma(\omega')$ for all $m \in \{0,\ldots,k-1\}$.  The inequalities $\sigma(\nu) \preceq \sigma^{m}(\omega')$ for $\omega'_{m+1} = 0$ and $\nu \prec \sigma^{m}(\omega')$ for $\omega'_{m+1} = 1$ follow from $\omega \in \Omega^{-}(\omega, \nu)$ and $\sigma^m(\omega) \prec \sigma^m(\omega')$.
To conclude the proof of the claim, it remains to show that  $\sigma^{m}(\omega') \preceq \omega'$ for all $m \in \{0,\ldots,k-1\}$ with $\omega'_{m+1} = 0$, as this implies that $\sigma^{m}(\omega') \preceq \sigma(\omega')$ when $\omega'_{m+1} = 1$. 
If $1 \le m < j$, then $\sigma^{m}(\omega) \preceq \omega$ and the minimality of $j$ imply that $(\omega_{m+1}, \dots, \omega_{n}) \prec (\omega_{1}, \dots, \omega_{n-m})$, hence $\sigma^{m}(\omega') \prec \omega'$. 
If $m \ge j$, then we have 
\begin{align*}
(\omega_{m+1}, \dots, \omega_{k}, \omega_{1}) \prec (\omega_{m+1}, \dots, \omega_{k}, 1) = (\omega_{m-j+1}, \dots, \omega_{k-j}, \omega_{k-j+1} ) \preceq (\omega_1, \dots, \omega_{k-m}, \omega_{k-m+1})
\end{align*}
when $\omega_{m+1} = 0$, which yields $\sigma^{m}(\omega') \prec \omega'$.
Therefore, we have $\omega' \in \Omega^{-}(\omega', \nu)$.
Since $n$ can be chosen arbitrarily large, we have that $\omega'$ can be made to be arbitrarily close to $\omega$ with respect to the metric $\mathscr{D}$.

Analogous to the construction of $\omega'$, one can build a periodic word $\nu' \preceq \nu$ arbitrarily close to $\nu$ with respect to the metric $\mathscr{D}$, such that $(\omega', \nu')$ is a pair of infinite words satisfying \Cref{defn:admissible} \eqref{enumi1:defn_admissible} and \eqref{enumi2:defn_admissible}.  By construction we have that $\Omega_{\beta,\alpha} \subset \Omega(\omega', \nu')$ and so
\begin{align*}
\ln(b) \coloneqq \lim_{k \to \infty} k^{-1} \ln \left( \lvert \Omega(\omega', \nu')\lvert_{k} \rvert \right) \ge \ln(\beta) > 0.
\end{align*}
Therefore, \Cref{defn:admissible} \eqref{enumi3:defn_admissible} is satisfied. 
By \Cref{rmk:subaddiitve}, if $\omega'\vert_{n} = \omega\vert_{n}$ and $\nu'\vert_{n} = \nu\vert_{n}$, then $\lvert \Omega_{\beta, \alpha}\lvert_{n} \rvert = \lvert \Omega(\omega', \nu')\lvert_{n} \rvert$ and thus $\ln(b) \leq n^{-1} \ln(\lvert \Omega_{\beta, \alpha}\lvert_{n} \rvert)$.  Hence, for a given $\epsilon > 0$, we can find an $N \in \mathbb{N}$ such that $0 \le b - \beta < \epsilon$, for all $n \geq N$.

Let $a$ be as in \eqref{eq:BSV:2014_alpha}, that is the solution of $\pi_{b,a}(\omega') = q$ with $q = (1 - a)/b$.  For $0 \le b-\beta < \epsilon$, we observe the following chain of inequalities.
\begin{align*}
\lvert \alpha - a \rvert
	&= \left\lvert 1 - \beta + \beta (\beta - 1 ) \sum_{i = 1}^{\infty} \omega_{i} \beta^{-i} - 1 + b - b (b - 1 )  \sum_{i = 1}^{\infty} \omega'_{i} b^{-i} \right\rvert\\
	&\leq b - \beta + \beta (\beta - 1) \sum_{i = n+1}^{\infty} \omega_{i} \beta^{-i} + b (b - 1) \sum_{i = n+1}^{\infty} \omega_{i} b^{-i} + \sum_{i = 1}^{n} \big\lvert  \beta(\beta - 1)\beta^{-i} - b(b - 1)b^{-i} \big\rvert\\
& \leq 2\, (b-\beta) + \beta^{-n+1} + b^{-n+1} + 2 \sum_{i = 1}^{n} (\beta^{-i} - b^{-i}) \\
& \leq 2\, \left(b-\beta + \beta^{-n+1} + \frac{1}{\beta-1} - \frac{1}{b-1} + \frac{1}{b^n(b-1)}\right)\\
&\leq 2\, \left(\epsilon + \frac{\epsilon}{(\beta-1)^2} + \frac{\beta^2-\beta+1}{\beta^n(\beta-1)}\right)\\
&\leq \frac{4\epsilon}{(\beta-1)^2} + \frac{6}{\beta^n(\beta-1)}
\end{align*}
This implies, given $\epsilon > 0$, there exists $n \in \mathbb{N}$ so that 
$\omega'$ and $\nu'$ are periodic, $\lvert b - \beta \rvert < \epsilon$ and $\lvert a - \alpha \rvert < \epsilon$.  If $\omega' = \tau^{-}_{b, a}(q)$ and $\nu' = \tau^{+}_{b, a}(q)$, then the intermediate $\beta$-shift $\Omega_{b, a}$ is a subshift of finite type.  Otherwise, an application of \Cref{lem:BSV:2014} yields the required result.

\section*{Acknowledgements}

B.\,Li, T. Sahlsten and T.\,Samuel are extremely grateful to the organisers of the \textsl{Research Program: Fractals and Dynamics} and the \textsl{Mittag-Leffler Institut} for their very kind hospitality during the writing of this article. T. Sahlsten was partially supported by the \textsl{Marie Sk{\l}odowska-Curie Individual Fellowship grant $\#\,655310$} and a travel grant from the \textsl{Finnish Academy of Science and Letters}.

\end{document}